\documentclass{article}    
\usepackage{amstext}
\usepackage{latexsym}

\begin{document}           

\begin{center}
{\huge When are There Infinitely Many Irreducible Elements in a Principal Ideal Domain?} \\ [.250in]
{\large Fabrizio Zanello}
\end{center}


\renewcommand{\baselinestretch}{1.50}
\setlength{\parskip}{10pt}

\def\qed{$\rlap{$\sqcap$}\sqcup$}

{\large


{\ }\\
\\
It has been a well-known fact since Euclid's time that there exist infinitely many rational primes. Two natural questions arise: In which other rings, sufficiently similar to the integers, are there infinitely many irreducible elements? Is there a unifying algebraic concept that characterizes such rings?\\
\indent The purpose of the present note is to place the fact concerning the infinity of primes into a more general context, one that also includes the interesting case of the factorial domains (unique factorization domains) of algebraic integers in a number field. We show that, if $A$ is a Principal Ideal Domain (PID, for short), then the fact that $A$ contains infinitely many (pairwise nonassociate) irreducible elements is equivalent to the property that every maximal ideal in the polynomial ring $A[x]$ has the same (maximal) height.

\indent We begin by recalling some basic definitions. (We assume that all rings are commutative with an identity element, denoted by 1). The {\it Jacobson radical} $J(R)$ of a ring $R$ is the intersection of all the maximal ideals of $R$, while the {\it nilradical} $\sqrt{0}$ of $R$ is the intersection of all the prime ideals of $R$. The latter can also be described as the set of all nilpotent elements of $R$ (see [1, Proposition 1.8, p. 5]). The {\it height} of a prime ideal $P$ in $R$ is the supremum of the lengths of the chains $P_0\subset P_1\subset ...\subset P_r=P$ of prime ideals of $R$ that end at $P$. The {\it Krull-dimension} $\dim R$ of $R$ is the supremum of the lengths of all the chains of prime ideals of $R$, or, equivalently, the supremum of the heights of all the prime ideals $P$ in $R$. For instance, a field has Krull-dimension 0, while a PID $A$ has Krull-dimension 1, since every prime ideal of $A$ different from $(0)$ has height 1. Finally, two elements $a$ and $b$ of $R$ are {\it associates} if there exists a unit $u$ in $R$ such that $a=ub$. Otherwise, we say that $a$ and $b$ are {\it nonassociates}.\\
\\
{\bf Lemma 1.} {\it Let $R$ be a domain. Then for the ring $R[x]$ it is the case that $J(R[x])=\sqrt{0}=(0).$}\\
\\
{\it Proof.} Of course $\sqrt{0}=(0)$, for $R[x]$ contains no nonzero nilpotent elements. Now let $f(x)$ be a member of $J(R[x])$. Then $xf(x)+1$ is a unit: otherwise $xf(x)+1$ would belong to some maximal ideal $M$ in $R[x]$, and since $f(x)$ lies in $M$, this would place 1 in $M$, a contradiction. Since $R[x]$ is a domain, this implies that $xf(x)+1$ is a constant polynomial, so $f=0$.{\ }{\ }\qed 

\indent Actually, the conclusion of Lemma 1 can be extended to an arbitrary ring $R$. It is not difficult to show (but the argument is longer) that in every $R[x]$ we have $J(R[x])=\sqrt{0}$ (see [1, Exercise 4, p. 11]).\\
\\
{\bf Theorem 2.} {\it Let $A$ be a PID. Then the following statements are equivalent:\\
(i) every maximal ideal of $A[x]$ has height 2;\\
(ii) $A$ has infinitely many pairwise nonassociate irreducible elements.}\\
\\
{\it Proof.} Let $M$ be a maximal ideal in $A[x]$ and set $P=M \cap A$. Then $P$ is a prime ideal of $A$. Suppose for the moment that $P=(0)$. We claim that $M$ has height 1 in this situation. Choose an element $g(x)$ in $M$ of least degree. Then $\deg g(x) >0$, and we can assume that $g(x)$ is irreducible, since $M$ is a prime ideal.

\indent Let $K$ be the quotient field of $A$. If $f(x)$ belongs to $M$, then in $K[x]$ we can write $$f(x)=g(x)q_1(x)+r_1(x),$$ where either $r_1(x)\equiv 0$ or $\deg r_1(x)<\deg g(x)$. Therefore we obtain $$f(x)=g(x){q(x)\over a}+{r(x)\over a},$$ or $$af(x)=g(x)q(x)+r(x),$$ with $q(x)$ and $r(x)$ in $A[x]$ and $a$ in $A-\lbrace 0\rbrace$. This implies that $r(x)$ lies in $M$, whence $r(x)\equiv 0$, because $\deg r(x)\geq \deg g(x)$ by the choice of $g(x)$. Therefore $af(x)\in (g(x))\subset M$. We infer that $f(x)$ belongs to $(g(x))$: $(g(x))$ is a prime ideal, so either $a$ or $f(x)$ belongs to it, but prime ideals in $A[x]$ cannot contain nonzero constants. Hence $M=(g(x))$, from which it follows that $M$ has height 1 in case $P=(0)$.\\
{\it Proof of (i) $\Rightarrow $(ii)}. Assume that every maximal ideal $M$ in $A[x]$ has height 2. From what we have just shown, it follows that every such $M$ must contain some nonzero constant, hence must include at least one irreducible element of $A$.

\indent Suppose now that there are only finitely many pairwise nonassociate irreducible elements of $A$ and that $p_1,p_2,...,p_s$ is a complete list of representatives of irreducible elements. Then $p_1p_2\cdot \cdot \cdot p_s$ belongs to every maximal ideal $M$ in $A[x]$, that is, to $J(A[x])$. Because $A$ is a domain, Lemma 1 asserts that $J(A[x])=(0)$, a contradiction. This proves (ii).\\
{\it Proof of (ii) $\Rightarrow $(i)}. Let $M$ again be an arbitrary maximal ideal of $A[x]$. We want to show that $M$ has height 2.

\indent If $P=M \cap A=(0)$, then we have already seen that $M=(g(x))$ for some irreducible $g(x)$ in $A[x]$. Since $M$ is a maximal ideal, for every irreducible $p$ in $A$ we must have $(g(x),p)=(1)$ (i.e., $A[x]/(g(x),p)=0$), and therefore $\overline{g(x)}$ must be a constant in $(A/(p))[x]$ for each such $p$. In other words, every coefficient of $g(x)$ other than the constant term is divisible by every irreducible element $p$ of $A$. This is a contradiction, since we are assuming that there are infinitely many pairwise nonassociate irreducible elements of $A$ and since, as a PID, $A$ is a factorial domain.

\indent Hence $P=M \cap A\neq (0)$, so $P=(p)$ for some irreducible element $p$ in $A$. Moreover $M$ is not principal: if it were, it would coincide with $(p)$, which is impossible because $A[x]/(p)\simeq (A/(p))[x]$ is not a field.\\
\indent Hence $M$ strictly contains the prime ideal $(p)$, and therefore its height is at least 2. But $$\dim A[x]=\dim A+1=2,$$ since $A$ is Noetherian (see [1, Exercise 7, p. 126]). Thus the height of $M$ is exactly 2, as required by (i). The proof of the theorem is complete.{\ }{\ }\qed \\
\\
{\bf ACKNOWLEDGEMENTS.} I wish to express my warm gratitude to Professor Ram Murty for an interesting discussion that led to the idea of writing this paper.\\
\\
\\
{\bf \huge References}\\
\\
$[1]$ {\ }M.F. Atiyah and I.G. Macdonald, {\it Introduction to Commutative Algebra}, Addison-Wesley, Reading, MA, 1969.\\
\\
\\
Fabrizio Zanello\\
Department of Mathematics and Statistics, Queen's University\\
Kingston, Ontario\\
K7L 3N6 Canada\\
fabrizio@mast.queensu.ca

}

\end{document}